\documentclass{article} 

\usepackage{amsfonts,amssymb,amsmath,amsthm}
\usepackage{lscape}
\usepackage[all, poly]{xy} 

\title{Centralizer algebras of the primitive unitary reflection group of order $96$}
\author{Masashi Kosuda and Manabu Oura} 

\theoremstyle{plain}
\newtheorem{thm}{Theorm}[section]
\newtheorem{cor}[thm]{Corollary}

\theoremstyle{definition}

\newcommand{\Z}{{\mathbb Z}}
\newcommand{\Comp}{\mathbb C}
\newcommand{\bolde}{\mbox{\boldmath $e$}}

\newcommand{\boldv}{\mbox{\boldmath $v$}}
\newcommand{\boldX}{\mbox{\boldmath $X$}}
\newcommand{\pprime}{\prime\prime}
\newcommand{\ppprime}{\prime\pprime}
\newcommand{\End}{\mbox{End}}

\newcommand{\FF}{{\mathbb F}}

\date{}

\begin{document}
\maketitle
\begin{abstract}
Among the unitary reflection groups,
the one on the title is singled out by its importance 
in, for example, coding theory and number theory.
In this paper we start with describing the irreducible representations of this group
and then examine
the semi-simple structure of the centralizer 
algebra in the tensor representation.
\end{abstract}

\renewcommand{\thefootnote}{\fnsymbol{footnote}}
\footnote[0]{
Keywords: Centralizer algebra, unitary reflection group, Bratteli diagram\newline
\hspace{0.4cm}
MSC2010:
Primary 05E10, Secondary 05E05 05E15 05E18.
\newline
\hspace{0.4cm}
running head: Centralizer algebra of $H_1$ in the tensor representation
}

\section{Introduction}
The group, which we denote by $H_1$, on the title of this paper
consists of $96$ matrices of size $2$ by $2$.
It is the unitary group generated by reflections (u.g.g.r.),
numbered as No.$8$ in Shephard-Todd~\cite{ST}.
This group,
as well as No.$9$ in the same list, has long been recognized.
The purpose of the present paper is to give a contribution to $H_1$
by decomposing the centralizer algebra of $H_1$
in the tensor representation into irreducible components.

We shall give an outline of the first statement in Abstract.
The group $H_1$ naturally acts on the polynomial ring 
$\Comp[x,y]$ of $2$ variables over the complex number field $\Comp$,
i.e.  
\[
Af(x,y)=f(ax+by,cx+dy), \ \ 
A=\begin{pmatrix}a & b \\ c & d \end{pmatrix}\in H_1
\]
for $f\in \Comp[x,y]$.
We consider the invariant ring 
\[
\Comp[x,y]^{H_1}
=\{f\in \Comp[x,y]:\
\ Af=f, \ \ \forall A\in H_1\}
\]
of $H_1$.
This ring has a rather simple structure.
It is generated by two algebraically 
independent homogeneous polynomials of degrees $8$ and $12$, and conversely
this nature characterizes the u.g.g.r.
Brou\'{e}-Enguehard~\cite{BE} found a map
connecting this invariant ring with number theory.
Take a homogeneous polynomial $f(x,y)$ of degree $n$
from the invariant ring.
Introducing theta constants
\[
\theta_{ab}(\tau)
=\sum_{ m\in \Z}
\exp 2\pi i
	\left[
		\frac{1}{2} \tau
		\left( m+\frac{a}{2} \right)^2
		+\left( m+\frac{a}{2} \right)
		\frac{b}{2}
	\right],
\]
we get a modular form
$f(\theta_{00}(2\tau),\theta_{10}(2\tau))$ 
of weight $n/2$ for $SL(2,\Z)$.
Moreover this map is an isomorphism from the
invariant ring of $H_1$ onto the ring of
modular forms for $SL(2,\Z)$.

Next we proceed to coding theory.
Let $\FF_2=\{0,1\}$ be the field of two elements
and $\FF_2^n$ the vector space of dimension $n$
over $\FF_2$
equipped with the usual inner product
$(u,v)=u_1v_1+\cdots + u_nv_n$.
The weight of a vector $u$ is the number of 
non-zero coordinates of $u$.
A code of length $n$ is by definition a linear subspace of $\FF_2^n$.
We impose two conditions on codes.
The first one is the self-duality which says that 
a code $C$ coincides with its dual code $C^{\perp}$,
that is, $C=C^{\perp}$
in which
\[
C^{\perp}=\{u \in \FF_2^n:\ 
(u,v)=0,\ \ \forall v\in C\}.
\]
The second one is the doubly-evenness which means
\[
wt(u)\equiv 0\pmod{4},\ \ \forall u\in C.
\]
These two notions give rise to the relation with invariant theory
via 
the weight enumerator 
\[
W_C(x,y)=\sum_{v\in C}
x^{n-wt(v)}y^{wt(v)}
\]
of a code $C$.
In fact, if $C$ is self-dual, we have
\[
W_C((x-y)/\sqrt{2},(x+y)/\sqrt{2})=W_C(x,y)
\]
and if $C$ is doubly even, we have
\[
W_C(x,iy)=W_C(x,y).
\]
We mention that a self-dual and doubly even code of length $n$ exists if and only if $n$ is a 
multiple of $8$.

Now we can state the connections among
all what we have mentioned.
Take a positive integer $n\equiv 0\pmod{8}$.
The weight enumerator of a self-dual doubly even code of length $n$
is an invariant of $H_1$
and 
\[
W_C(\theta_{00}(2\tau),\theta_{10}(2\tau))
\]
is a modular form of weight $n/2$ for $SL(2,\Z)$.
Gleason~\cite{Gl} showed that
the invariants of degree $n$ can be
spanned by the weight enumerators of self-dual 
doubly even codes of length $n$.
Finally any modular form of weight $n/2$ can be obtained from 
the weight enumerator of self-dual doubly even codes of length $k$.
The whole theory with more general results could be found in \cite{Ru1}, \cite{Ru2}
from which our notation $H_1$ comes.

Besides the importance of $H_1$,
the motivation of this paper could be found 
in Brauer~\cite{Br}, Weyl~\cite{We}.
One of the main ingredients there is the {\it commutator algebra}
where invariant theory comes into play.
We follow Weyl.
Given any group of linear transformations in an $n$-dimensional space.
Take covariant vectors $y^{(1)},\ldots,y^{(f)}$ 
and contravariant vectors $\xi^{(1)},\ldots,\xi^{(f)}$.
A linear transformation acts on covariant vectors {\it cogrediently}
and on contravariant vectors {\it contragradiently}.
Then
the matrices $\| b(i_1\cdots i_f;\ k_1\cdots k_f)\|$ in the tensor space
obtained from the invariants
\[
\sum_{i;k}
b(i_1\cdots i_f;\ k_1\cdots k_f)\xi_{i_1}^{(1)}
\cdots \xi_{i_f}^{(f)}y_{k_1}^{(1)}\cdots y_{k_f}^{(f)}
\]
form the commutator algebra of $H_1$ in the tensor representation.
The problem here is to decompose this algebra into simple parts.
It is quite natural to apply this philosophy to our group $H_1$
as we will in this paper ({\it cf}. \cite{Ba}).

\section{Irreducible representations of $H_1$}
In this section we determine the irreducible representations of $H_1$
which yields the character table.
At the end of this section we discuss invariant theory of $H_1$ under the 
irreducible representations.

The unitary reflection group $H_1$
is a finite group in $U_2$
generated by the following matrices $T$ and $D$:
\[
	T = \frac{1+i}{2}
	\begin{pmatrix}
		1 & 1\\
		1 & -1
	\end{pmatrix}
 = \frac{1}{\sqrt{2}}
	\begin{pmatrix}
		\epsilon & \epsilon\\
		\epsilon & \epsilon^5
	\end{pmatrix}, \quad
	D =
	\begin{pmatrix}
		1 & 0\\
		0 & i
	\end{pmatrix}
 .
\]
Here $\epsilon = \exp(2\pi i/8)$.
It is known that
the group size of $H_1$ is 96
and it has 16 conjugacy classes
$\mathfrak{C}_{1}, \ldots, \mathfrak{C}_{16}$.
Each conjugacy class has the following representative:
\begin{eqnarray*}
&\mathfrak{C}_{1} \ni 1
=\begin{pmatrix}1&0\\0&1\end{pmatrix},\ %
\mathfrak{C}_{2} \ni T
= \frac{1}{\sqrt{2}}
	\begin{pmatrix}
		\epsilon & \epsilon\\
		\epsilon & \epsilon^5
	\end{pmatrix},\ %
\mathfrak{C}_{3} \ni T^2
=\begin{pmatrix}i&0\\0&i\end{pmatrix},
&\\
&\mathfrak{C}_{4} \ni T^3
= \frac{1}{\sqrt{2}}
	\begin{pmatrix}
		\epsilon^3 & \epsilon^3\\
		\epsilon^3 & \epsilon^7
	\end{pmatrix},\ %
\mathfrak{C}_{5} \ni T^4
=\begin{pmatrix}-1&0\\0&-1\end{pmatrix},\ %
\mathfrak{C}_{6} \ni T^6
=\begin{pmatrix}-i&0\\0&-i\end{pmatrix},
&\\
&\mathfrak{C}_{7} \ni D
=\begin{pmatrix}1&0\\0&i\end{pmatrix},\ %
\mathfrak{C}_{8} \ni DT
= \frac{1}{\sqrt{2}}
	\begin{pmatrix}
		\epsilon & \epsilon\\
		\epsilon^3 & \epsilon^7
	\end{pmatrix},\ %
\mathfrak{C}_{9} \ni DT^2
=\begin{pmatrix}i&0\\0&-1\end{pmatrix},
&\\ 
&\mathfrak{C}_{10} \ni DT^3
= \frac{1}{\sqrt{2}}
	\begin{pmatrix}
		\epsilon^3 & \epsilon^3\\
		\epsilon^5 & \epsilon
	\end{pmatrix},\ %
\mathfrak{C}_{11} \ni DT^4
=\begin{pmatrix}-1&0\\0&-i\end{pmatrix},
&\\
&\mathfrak{C}_{12} \ni DT^5
=\frac{1}{\sqrt{2}}
	\begin{pmatrix}
		\epsilon^5 & \epsilon^5\\
		\epsilon^7 & \epsilon^3
	\end{pmatrix},\ %
\mathfrak{C}_{13} \ni DT^6
=\begin{pmatrix}-i&0\\0&1\end{pmatrix},
&\\
&\mathfrak{C}_{14} \ni DT^7
= \frac{1}{\sqrt{2}}
	\begin{pmatrix}
		\epsilon^7 & \epsilon^7\\
		\epsilon & \epsilon^5
	\end{pmatrix},\ %
\mathfrak{C}_{15} \ni D^2
=\begin{pmatrix}1&0\\0&-1\end{pmatrix},\ 
\mathfrak{C}_{16} \ni D^2T^2
=\begin{pmatrix}i&0\\0&-i\end{pmatrix}.
&
\end{eqnarray*}

Since the number of conjugacy classes
and that of the non-isomorphic irreducible representations
coincide,
there exist 16 classes of the irreducible representations of $H_1$.
In the following, we  construct all of them one by one.

First we note that any group has the trivial representation
which maps each element of the group to 1.
We denote that of $H_1$ by $(\rho_1, V_1)$.
The determinant which maps $T$ and $D$ to $-i$ and $i$ respectively
also gives a one-dimensional irreducible representation.
We call it $(\rho_3, V_3)$.
The tensor product $\rho_3^{\otimes 2}$ also
gives a one-dimensional irreducible representation,
which maps both $T$ and $D$ to $-1$.
We name it $(\rho_2, V_2)$.
Also $\rho_2\otimes\rho_3$ defines a one-dimensional representation.
We name it $(\rho_4, V_4)$.

Next we consider two-dimensional representations.
The natural representation $(\rho_{10}, V_{10})$
which maps $T$ and $D$ to the defining matrices above
is irreducible,
since neither of one-dimensional $D$-invariant subspaces
are $T$-invariant. 
Taking tensor products with the one-dimensional representations above
and the natural representation,
we have further 3 two-dimensional irreducible representations,
$\rho_7 = \rho_3\otimes\rho_{10}$,
$\rho_8 = \rho_2\otimes\rho_{10}$ and $\rho_9 = \rho_4\otimes\rho_{10}$.
There are 2 more two-dimensional irreducible representations
which we will deal with later.

As a subrepresentation of $\rho_{10}\otimes\rho_{10}$,
we have a three-dimensional irreducible representation.
Let $\langle \bolde_1, \bolde_2 \rangle$ be a basis of $V_{10}$
which gives the natural representation.
Then $\langle \bolde_1\otimes\bolde_1,
	\bolde_1\otimes\bolde_2,
	\bolde_2\otimes\bolde_1,
	\bolde_2\otimes\bolde_2\rangle$
gives a basis for the tensor representation $\rho_{10}^{\otimes 2}$.
With respect to this basis, the representation matrices of $T$ and $D$ are
\[
	\rho_{10}^{\otimes 2}(T)
	= \frac{i}{2}
	\begin{pmatrix}
		1& 1& 1& 1\\
		1&-1& 1&-1\\
		1& 1&-1&-1\\
		1&-1&-1& 1
	\end{pmatrix}\ \mbox{ and }\ 
	\rho_{10}^{\otimes 2}(D)
	= \mbox{diag}(1, i, i, -1).
\]
If we put $\bolde'_1 = \bolde_1\otimes\bolde_1$,
$\bolde'_2 = \bolde_1\otimes\bolde_2+\bolde_2\otimes\bolde_1$
and $\bolde'_3 = \bolde_2\otimes\bolde_2$,
then $\langle \bolde'_1, \bolde'_2, \bolde'_3\rangle$
is obviously a $D$-invariant subspace.
It is easy to check that it is also $T$-invariant.
Hence it gives a three-dimensional representation.
We name it ($\rho_{13}$, $V_{13}$).
The representation matrices with respect to this basis are
\[
	\rho_{13}(T)
	= \frac{i}{2}
	\begin{pmatrix}
		1& 2& 1\\
		1& 0&-1\\
		1&-2& 1
	\end{pmatrix}\ \mbox{and}\ 
	\rho_{13}(D)
	= \mbox{diag}(1, i, -1).
\]
Since each one-dimensional $D$-invariant subspace of $V_{13}$
is not $T$-invariant,
the representation ($\rho_{13}$, $V_{13}$) is irreducible.
Similarly to the previous case,
we have further 3 three-dimensional irreducible representations,
$\rho_{11} = \rho_3\otimes\rho_{13}$,
$\rho_{12} = \rho_4\otimes\rho_{13}$ and
$\rho_{14} = \rho_2\otimes\rho_{13}$.

Next we look for a four-dimensional irreducible representation
in $(\rho_{10}\otimes\rho_{13}, V_{10}\otimes V_{13})$.
Let $\langle \bolde_i\otimes\bolde'_j\ |\ i=1,2,\ j=1,2,3\rangle$
be a basis of $V_{10}\otimes V_{13}$ (lexicographical order).
Then we have the following representation matrices of $T$ and $D$:
\begin{align*}
	\rho_{10}\otimes\rho_{13}(T)
	&= \frac{-1+i}{4}
	\begin{pmatrix}
		1& 2& 1& 1& 2& 1\\
		1& 0&-1& 1& 0&-1\\
		1&-2& 1& 1&-2& 1\\
		1& 2& 1&-1&-2&-1\\
		1& 0&-1&-1& 0& 1\\
		1&-2& 1&-1& 2&-1
	\end{pmatrix},\\
	\rho_{10}\otimes\rho_{13}(D)
	&= \mbox{diag}(1, i, -1, i,-1, -i).
\end{align*}
If we put $\bolde^{\pprime}_1 = \bolde_1\otimes\bolde'_1$,
$\bolde^{\pprime}_2 = \bolde_1\otimes\bolde'_2+\bolde_2\otimes\bolde'_1$,
$\bolde^{\pprime}_3 = \bolde_1\otimes\bolde'_3+\bolde_2\otimes\bolde'_2$,
and $\bolde^{\pprime}_4 = \bolde_2\otimes\bolde'_3$,
according to the eigen values of $\rho_{10}\otimes\rho_{13}(D)$,
then $\langle \bolde^{\pprime}_k\ |\ k=1,2,3,4 \rangle$
is obviously a $D$-invariant subspace.
It is also easy to check that it is $T$-invariant.
Hence it gives a four-dimensional representation.
We name it ($\rho_{15}$, $V_{15}$).
The representation matrices with respect to this basis are
\[
	\rho_{15}(T)
	= \frac{-1+i}{4}
	\begin{pmatrix}
		1& 3& 3& 1\\
		1& 1&-1&-1\\
		1&-1&-1& 1\\
		1&-3& 3&-1
	\end{pmatrix}\ \mbox{and}\ 
	\rho_{13}(D)
	= \mbox{diag}(1, i, -1, -i).
\]
As we saw in the previous case,
none of one-dimensional $D$-invariant subspaces of $V_{15}$
is $T$-invariant.
Now consider two-dimensional $D$-invariant subspaces.
Since all eigen spaces of $\rho_{15}(D)$ are one-dimensional,
we find that a two-dimensional $D$-invariant subspace
is of the form
$\langle \bolde^{\pprime}_i, \bolde^{\pprime}_j\rangle (i\neq j)$.
Let $W$ be $\langle \bolde^{\pprime}_1, \bolde^{\pprime}_2\rangle$
and take a non-zero vector
$\boldv = a\bolde^{\pprime}_1+b\bolde^{\pprime}_2$ from $W$.
Then we have
\[
	\rho_{15}(T)\boldv
	= \frac{-1+i}{4}
	\left[
		(a+3b)\bolde^{\pprime}_1
		+ (a+b)\bolde^{\pprime}_2
		+ (a-b)\bolde^{\pprime}_3
		+ (a-3b)\bolde^{\pprime}_4
	\right].
\]
In order that $\rho_{15}(T)\boldv\in W$,
it must hold that $a=b=0$.
This contradicts the assumption that $\boldv$ is non-zero vector.
Hence we find that $W$ is not $T$-invariant.
Similar arguments hold
for all two-dimensional $D$-invariant subspaces
$\{\langle \bolde^{\pprime}_i,
\bolde^{\pprime}_j\rangle\}_{1\leq i < j\leq 4}$.
This implies there is no two-dimensional subrepresentation in $V_{15}$.
Hence we find that ($\rho_{15}$, $V_{15}$) is irreducible.
Similarly to the previous case,
we have further four-dimensional irreducible representations,
$\rho_2\otimes\rho_{15}$, $\rho_3\otimes\rho_{15}$, $\rho_4\otimes\rho_{15}$. 
The first one, however, coincides with $\rho_{15}$
and the second one and the third one are equivalent.
Hence we have 2 four-dimensional irreducible representations,
$\rho_{15}$ and $\rho_{16} = \rho_3\otimes\rho_{15}$.

Finally we look for the remaining irreducible representations
in $(\rho_{10}\otimes\rho_{15}, V_{10}\otimes V_{15})$.
Let $\langle \bolde_i\otimes\bolde^{\pprime}_j\ |\ i=1,2, j=1,2,3,4\rangle$
be a basis of $V_{10}\otimes V_{15}$ (lexicographical order).
Then we have the following representation matrices of $T$ and $D$:
\begin{align*}
	\rho_{10}\otimes\rho_{15}(T)
	&= \frac{-1}{4}
	\begin{pmatrix}
		1& 3& 3& 1& 1& 3& 3& 1\\
		1& 1&-1&-1& 1& 1&-1&-1\\
		1&-1&-1& 1& 1&-1&-1& 1\\
		1&-3& 3&-1& 1&-3& 3&-1\\
		1& 3& 3& 1&-1&-3&-3&-1\\
		1& 1&-1&-1&-1&-1& 1& 1\\
		1&-1&-1& 1&-1& 1& 1&-1\\
		1&-3& 3&-1&-1& 3&-3& 1
	\end{pmatrix},\\
	\rho_{10}\otimes\rho_{15}(D)
	&= \mbox{diag}(1, i, -1,-i, i, -1,-i, 1).
\end{align*}
If we put $\bolde^{\ppprime}_1
= \bolde_1\otimes\bolde^{\pprime}_1+\bolde_2\otimes\bolde^{\pprime}_4$,
and $\bolde^{\ppprime}_2
= \bolde_1\otimes\bolde^{\pprime}_3+\bolde_2\otimes\bolde^{\pprime}_2$,
then $\langle \bolde^{\ppprime}_1, \bolde^{\ppprime}_2\rangle$
is $T$- and $D$-invariant subspace.
We name it ($\rho_{5}$, $V_{5}$).
The representation matrices with respect to this basis are
\[
	\rho_{5}(T)
	= \frac{-1}{2}
	\begin{pmatrix}
		1& 1\\
		3&-1
	\end{pmatrix}\ \mbox{and}\ 
	\rho_{5}(D)
	= \mbox{diag}(1, -1).
\]
Similarly to the previous ones,
we can check that this representation is irreducible.
Further $\rho_6 = \rho_3\otimes\rho_5$
also defines a two-dimensional irreducible representation.

So far, we have got 16 irreducible representations.
Since $H_1$ has 16 conjugacy classes,
$\{(\rho_i, V_i)\}_{i=1}^{16}$
are a complete representatives of all irreducible representations
of $H_1$.
Accordingly,
the character table of $H_1$ is also derived.

\begin{landscape}
\begin{center}
\begin{tabular}{c|cccccccccccccccc}
\hline
$H_1$&	$\mathfrak{C}_1$&$\mathfrak{C}_2$&$\mathfrak{C}_3$&$\mathfrak{C}_4$&
		$\mathfrak{C}_5$&$\mathfrak{C}_6$&$\mathfrak{C}_7$&$\mathfrak{C}_8$&
	$\mathfrak{C}_9$&   $\mathfrak{C}_{10}$&$\mathfrak{C}_{11}$&
	$\mathfrak{C}_{12}$&$\mathfrak{C}_{13}$&$\mathfrak{C}_{14}$&
	$\mathfrak{C}_{15}$&$\mathfrak{C}_{16}$\\
     &       $1$&     $T$&   $T^2$&   $T^3$&
	   $T^4$&   $T^6$&     $D$&    $DT$&
	  $DT^2$&  $DT^3$&  $DT^4$&  $DT^5$&
	  $DT^6$&  $DT^7$&   $D^2$&$D^2T^2$\\
order	& 1& 8& 4& 8&
	  2& 4& 4& 6&
	  4&12& 4& 3&
	  4&12& 2& 4\\
\hline
$\chi_1$&1&1&1&1&
	 1&1&1&1&
	 1&1&1&1&
	 1&1&1&1\\
$\chi_2$&  1&$-1$&   1&$-1$&
	   1&   1&$-1$&   1&
	$-1$&   1&$-1$&   1&
	$-1$&   1&   1&   1\\
$\chi_3$&  1&$-i$&$-1$& $i$&
	   1&$-1$& $i$&   1&
	$-i$&$-1$& $i$&   1&
	$-i$&$-1$&$-1$&   1\\
$\chi_4$&  1& $i$&$-1$&$-i$&
	   1&$-1$&$-i$&   1&
	 $i$&$-1$&$-i$&   1&
	 $i$&$-1$&$-1$&   1\\
$\chi_5$&  2&   0&   2&   0&
	   2&   2&   0&$-1$&
	   0&$-1$&   0&$-1$&
	   0&$-1$&   2&   2\\
$\chi_6$&  2&   0&$-2$&   0&
	   2&$-2$&   0&$-1$&
	   0&   1&   0&$-1$&
	   0&   1&$-2$&   2\\
$\chi_7$&     2&   0& $-2i$&   0&
	   $-2$&$2i$&$-1+i$&   1&
	  $1+i$&$-i$& $1-i$&$-1$&
	 $-1-i$& $i$&     0&   0\\
$\chi_8$&     2&    0&  $2i$&   0&
	   $-2$&$-2i$&$-1-i$&   1&
	  $1-i$&  $i$& $1+i$&$-1$&
	 $-1+i$& $-i$&     0&   0\\
$\chi_9$&    2&   0& $-2i$&   0&
	  $-2$&$2i$& $1-i$&   1&
	$-1-i$&$-i$&$-1+i$&$-1$&
	 $1+i$& $i$&     0&   0\\
$\chi_{10}$& 2&    0&  $2i$&   0&
	  $-2$&$-2i$& $1+i$&   1&
	$-1+i$&  $i$&$-1-i$&$-1$&
	 $1-i$& $-i$&     0&   0\\
$\chi_{11}$& 3&1&   3&  1&
	     3&3&$-1$&  0&
	  $-1$&0&$-1$&  0&
	  $-1$&0&$-1$&$-1$\\
$\chi_{12}$& 3&$-1$&   3&$-1$&
	     3&   3&   1&   0&
	     1&   0&   1&   0&
	     1&   0&$-1$&$-1$\\
$\chi_{13}$&    3& $i$&$-3$&$-i$&
	        3&$-3$& $i$&   0&
	     $-i$&   0& $i$&   0&
	     $-i$&   0&1&$-1$\\
$\chi_{14}$&   3&$-i$& $-3$& $i$&
	       3&$-3$& $-i$&   0&
	     $i$&   0& $-i$&   0&
	     $i$&   0&    1&$-1$\\
$\chi_{15}$&   4&   0&$-4i$&   0&
	    $-4$&$4i$&    0&$-1$&
	       0& $i$&    0&   1&
	       0&$-i$&    0&   0\\
$\chi_{16}$&   4&    0&$4i$&   0&
	    $-4$&$-4i$&   0&$-1$&
	       0& $-i$&   0&   1&
	       0&  $i$&   0&   0\\
\end{tabular}
\end{center}
\end{landscape}

We conclude this section with adding a few words on invariant theory of $H_1$ under irreducible representations
({\it cf}. \cite{Di}).
Let $\rho$ be one of the $d$-dimensional irreducible representation of $H_1$.
Then $\rho (H_1)$ acts naturally on the polynomial ring of $d$ variables.
We denote the invariant ring under this action by $\Comp [\rho]^{H_1}$.
The orders of $\rho_i (H_1)$ are
\[
\underbrace{1,2,4,4}_{\dim 1},
\underbrace{6,12,96,96,96,96}_{\dim 2},
\underbrace{24,24,48,48}_{\dim 3},
\underbrace{96,96}_{\dim 4}.
\]
The dimension $1$ case aside, the invariant rings
\[
\Comp [\rho_5]^{H_1},
\Comp [\rho_i]^{H_1}\ (i=7,8,9,10),
\Comp [\rho_{12}]^{H_1}
\]
are weighted polynomial rings.
In the sense of \cite{ST},
all $\rho_i(H_1)\  (i=7,8,9,10)$ are equivalent each other,
and $\rho_{15}(H_1)$ to $\rho_{16}(H_1)$.
We already know the ring $\Comp [\rho_7]^{H_1}$.
The ring $\Comp [\rho_5]^{H_1}$ can be generated by the polynomials of degrees $2$ and $3$,
and the ring $\Comp [ \rho_{15}]^{H_1}$ by those of degrees $2,3$ and $4$.
If we look at the degrees, we can find that
$\rho_5(H_1)$ is equivalent to $G(3,3,2)$
and $\rho_{12}(H_1)$ to $G(2,2,3)$.
The other cases up to dimension $3$ are modules of rank $2$ over the polynomial rings.
The $\rho_{15}$ case has a somewhat complicated structure.
The ring $\Comp[\rho_{15}]^{H_1}$ is a module of rank $32$ over the polynomial ring.
We note that calculations here were done with Magma~\cite{Bo}.

\section{Decomposition of tensor representations}
In the previous section,
we have found complete representatives of
all irreducible representations.
In this section, we see how tensor powers of $\rho_{10}$
are decomposed into irreducible ones.

We begin with the general theory (see for example Curtis-Reiner\cite{CR}).
Let $\chi_1,\ldots, \chi_s$ be the
set of all irreducible characters of a finite group $G$.
For any (not necessarily irreducible)
representation $(\rho, V)$ of $G$,
let $\chi$ be its character.
Then $\chi$ can be uniquely expressed a sum of irreducible characters:
\[
	\chi = m_1\chi_1+\cdots + m_s\chi_s.
\]
Now suppose that $\chi$
has its character values ($k_1, \ldots, k_s$)
on the conjugacy classes ($\mathfrak{C}_1, \ldots, \mathfrak{C}_s$).
Then we get
\begin{align*}
	(k_1,\ldots, k_s)
	&= (\chi(\mathfrak{C}_1), \ldots, \chi(\mathfrak{C}_s))\\
	&= m_1(\chi_1(\mathfrak{C}_1), \ldots, \chi_s(\mathfrak{C}_s))
	+\cdots
	+ m_s(\chi_s(\mathfrak{C}_1), \ldots, \chi_s(\mathfrak{C}_s))\\
	&= (m_1, \ldots, m_s)
		\begin{pmatrix}
			\chi_1(\mathfrak{C}_1)&\cdots&\chi_1(\mathfrak{C}_s)\\
			\vdots	&\ddots&	\vdots\\
			\chi_s(\mathfrak{C}_1)&\cdots&\chi_s(\mathfrak{C}_s)
		\end{pmatrix}.
\end{align*}
If we let $\boldX$ denote the matrix of the character table,
then the above relation is simply written as
\begin{equation}
	(k_1, \ldots, k_s) = (m_1, \ldots, m_s)\boldX.
\end{equation}
By the linear independence of the irreducible characters,
$\boldX$ is non-singular.
Hence we have
\begin{equation}\label{eq:MulIrr}
	(m_1, \ldots, m_s) = (k_1, \ldots, k_s)\boldX^{-1}.
\end{equation}

In order to examine the structure of the centralizer algebra
of the tensor representation,
it is useful to investigate how the tensor product
of the natural and an irreducible representation
is decomposed into the irreducible ones.
In the following, we go back to our case and decompose
$\rho_{10}\otimes\rho_i$ ($i = 1, 2, \ldots, 16$)
one by one.

By the argument and/or the character table in the previous section,
we already have the following:
\begin{align*}
	\chi_{10}\cdot\chi_1 &= \chi_{10},\\
	\chi_{10}\cdot\chi_2 &= \chi_{8},\\
	\chi_{10}\cdot\chi_3 &= \chi_{7},\\
	\chi_{10}\cdot\chi_4 &= \chi_{9}.
\end{align*}

Further, we can directly read the following from the character table:
\begin{align*}
	\chi_{10}\cdot\chi_5 &= \chi_{16},\\
	\chi_{10}\cdot\chi_6 &= \chi_{15}.
\end{align*}

Next, consider
$\chi_{10}\cdot\chi_7(\mathfrak{C}_{1}, \ldots, \mathfrak{C}_{16})$.
Again from the character table,
we have
\begin{align*}
\lefteqn{
	\chi_{10}\cdot\chi_7(\mathfrak{C}_{1}, \ldots, \mathfrak{C}_{16})
}\\
&= (4, 0, 4, 0, 4, 4,-2, 1,-2, 1, -2, 1, -2, 1, 0, 0).
\end{align*}
Using the identity \eqref{eq:MulIrr},
we have
\begin{align*}
\lefteqn{(0, 1, 0, 0, 0, 0, 0, 0, 0, 0, 1, 0, 0, 0, 0, 0)}\\
&= \chi_{10}\cdot
\chi_7(\mathfrak{C}_{1}, \ldots, \mathfrak{C}_{16})\boldX^{-1}.
\end{align*}
This means
\[
	\chi_{10}\cdot\chi_7 = \chi_2 + \chi_{11}.
\]
In a similar way we have
\begin{align*}
\chi\cdot\chi_8 &= \chi_4+\chi_{14}, \\
\chi\cdot\chi_9 &= \chi_1+\chi_{12}, \\
\chi\cdot\chi_{10} &= \chi_3 +  \chi_{13}, \\
\chi\cdot\chi_{11} &= \chi_8 +  \chi_{16}, \\
\chi\cdot\chi_{12} &= \chi_{10} +  \chi_{16}, \\
\chi\cdot\chi_{13} &= \chi_7 +  \chi_{15}, \\
\chi\cdot\chi_{14} &= \chi_9 +  \chi_{15}, \\
\chi\cdot\chi_{15} &= \chi_5 +  \chi_{11} + \chi_{12}, \\
\chi\cdot\chi_{16} &= \chi_6 +  \chi_{13} + \chi_{14}. \\
\end{align*}
By the above calculation,
we obtain the Hasse diagram of the decomposition of
$\rho_{10}^{\otimes k}$ into irreducible ones.
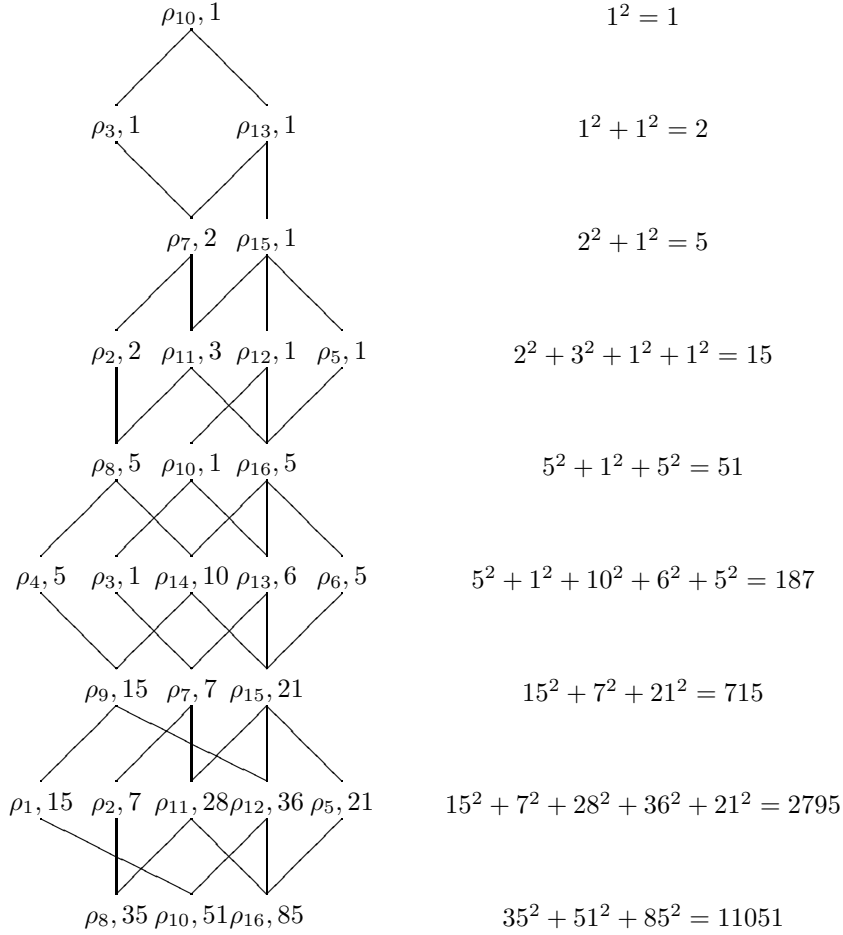
\[
\begin{xy}
(20,0)="T1P10",
(20,-3)*{\rho_{10},1},
(80,-3)*{1^2 = 1}, 
(20,-5)="S1P10",
(10,-15)="T2P3", (30,-15)="T2P13",
(10,-18)*{\rho_{3},1},(30,-18)*{\rho_{13},1},
(80,-18)*{1^2+1^2 = 2}, 
(10,-20)="S2P3", (30,-20)="S2P13",
(20,-30)="T3P7",(30,-30)="T3P15",
(20,-33)*{\rho_7,2}, (30,-33)*{\rho_{15},1},
(80,-33)*{2^2+1^2 = 5}, 
(20,-35)="S3P7",(30,-35)="S3P15",
(10,-45)="T4P2",(20,-45)="T4P11",(30,-45)="T4P12",(40,-45)="T4P5",(10,-48)*{\rho_2,2},(20,-48)*{\rho_{11},3},(30,-48)*{\rho_{12},1},(40,-48)*{\rho_5,1},
(80,-48)*{2^2+3^2+1^2+1^2 = 15}, 
(10,-50)="S4P2",(20,-50)="S4P11",(30,-50)="S4P12",(40,-50)="S4P5",(10,-60)="T5P8",(20,-60)="T5P10",(30,-60)="T5P16",
(10,-63)*{\rho_8,5},(20,-63)*{\rho_{10},1},(30,-63)*{\rho_{16},5},
(80,-63)*{5^2+1^2+5^2 = 51}, 
(10,-65)="S5P8",(20,-65)="S5P10",(30,-65)="S5P16", 
(0,-75)="T6P4",(10,-75)="T6P3",(20,-75)="T6P14",(30,-75)="T6P13",(40,-75)="T6P6",
(0,-78)*{\rho_4,5},(10,-78)*{\rho_{3},1},(20,-78)*{\rho_{14},10},
(30,-78)*{\rho_{13},6},(40,-78)*{\rho_{6},5},
(80,-78)*{5^2+1^2+10^2+6^2+5^2 = 187}, 
(0,-80)="S6P4",(10,-80)="S6P3",(20,-80)="S6P14",(30,-80)="S6P13",(40,-80)="S6P6", 
(10,-90)="T7P9",(20,-90)="T7P7",(30,-90)="T7P15",
(10,-93)*{\rho_{9},15},(20,-93)*{\rho_{7},7},
(30,-93)*{\rho_{15},21},
(80,-93)*{15^2+7^2+21^2 = 715}, 
(10,-95)="S7P9",(20,-95)="S7P7",(30,-95)="S7P15", 
(0,-105)="T8P1",(10,-105)="T8P2",(20,-105)="T8P11",(30,-105)="T8P12",(40,-105)="T8P5",
(0,-108)*{\rho_{1},15},(10,-108)*{\rho_{2},7},
(20,-108)*{\rho_{11},28},(30,-108)*{\rho_{12},36},
(40,-108)*{\rho_{5},21},
(80,-108)*{15^2+7^2+28^2+36^2+21^2 = 2795}, 
(0,-110)="S8P1",(10,-110)="S8P2",(20,-110)="S8P11",(30,-110)="S8P12",(40,-110)="S8P5", 
(10,-120)="T9P8",(20,-120)="T9P10",(30,-120)="T9P16",
(10,-123)*{\rho_{8},35},(20,-123)*{\rho_{10},51},
(30,-123)*{\rho_{16},85},
(80,-123)*{35^2+51^2+85^2 = 11051}, 
{\ar @{{.}-{.}}"S1P10";"T2P3"},{\ar @{{.}-{.}}"S1P10";"T2P13"},
{\ar @{{.}-{.}}"S2P3";"T3P7"},{\ar @{{.}-{.}}"S2P13";"T3P7"}, 
{\ar @{{.}-{.}}"S2P13";"T3P15"},
{\ar @{{.}-{.}}"S3P7";"T4P2"},{\ar @{{.}-{.}}"S3P7";"T4P11"}, 
{\ar @{{.}-{.}}"S3P15";"T4P11"},{\ar @{{.}-{.}}"S3P15";"T4P12"},
{\ar @{{.}-{.}}"S3P15";"T4P5"},
{\ar @{{.}-{.}}"S4P2";"T5P8"}, 
{\ar @{{.}-{.}}"S4P11";"T5P8"},{\ar @{{.}-{.}}"S4P11";"T5P16"},
{\ar @{{.}-{.}}"S4P12";"T5P10"},{\ar @{{.}-{.}}"S4P12";"T5P16"},
{\ar @{{.}-{.}}"S4P5";"T5P16"},
{\ar @{{.}-{.}}"S5P8";"T6P4"},{\ar @{{.}-{.}}"S5P8";"T6P14"}, 
{\ar @{{.}-{.}}"S5P10";"T6P3"},{\ar @{{.}-{.}}"S5P10";"T6P13"},
{\ar @{{.}-{.}}"S5P16";"T6P14"},{\ar @{{.}-{.}}"S5P16";"T6P13"},
{\ar @{{.}-{.}}"S5P16";"T6P6"},
{\ar @{{.}-{.}}"S6P4";"T7P9"}, 
{\ar @{{.}-{.}}"S6P3";"T7P7"},
{\ar @{{.}-{.}}"S6P14";"T7P9"},{\ar @{{.}-{.}}"S6P14";"T7P15"},
{\ar @{{.}-{.}}"S6P13";"T7P7"},{\ar @{{.}-{.}}"S6P13";"T7P15"},
{\ar @{{.}-{.}}"S6P6";"T7P15"},
{\ar @{{.}-{.}}"S7P9";"T8P1"},{\ar @{{.}-{.}}"S7P9";"T8P12"}, 
{\ar @{{.}-{.}}"S7P7";"T8P2"},{\ar @{{.}-{.}}"S7P7";"T8P11"},
{\ar @{{.}-{.}}"S7P15";"T8P11"},{\ar @{{.}-{.}}"S7P15";"T8P12"},
{\ar @{{.}-{.}}"S7P15";"T8P5"},
{\ar @{{.}-{.}}"S8P1";"T9P10"}, 
{\ar @{{.}-{.}}"S8P2";"T9P8"},
{\ar @{{.}-{.}}"S8P11";"T9P8"},{\ar @{{.}-{.}}"S8P11";"T9P16"},
{\ar @{{.}-{.}}"S8P12";"T9P10"},{\ar @{{.}-{.}}"S8P12";"T9P16"},
{\ar @{{.}-{.}}"S8P5";"T9P16"},
\end{xy}
\]
From this diagram we can read off the multiplicity of
each irreducible representation $(\rho_i, V_i)$ in $V_{10}^{\otimes k}$
by counting the number of paths from the top vertex
indexed by $\rho_{10}$ in the 1-st row
to the corresponding vertex in the $k$-th row.
We put the multiplicity on the right side of each irreducible
representation.
Further, we calculated the square sums of the multiplicities
on the each row.

Let ${\cal A}_k = \End_{H_1}(V_{10}^{\otimes k})$
be the centralizer algebra of $H_1$ in $V_{10}^{\otimes k}$,
where $H_1$ acts on $V_{10}$ diagonally.
By the Schur-Weyl reciprocity~\cite{Sa,We},
this diagram is the Bratteli diagram of
the algebra sequence
\[
	\Comp ={\cal A}_1\subset{\cal A}_2\subset{\cal A}_3\subset \cdots .
\]
(For the Bratteli diagram,
see for example Goodman-de la Harpe-Jones~\cite{GHJ}, \S 2.3.)
Accordingly,
the square sum of the multiplicities on the $k$-th row
is the dimension of $\mbox{End}_{H_1}(V_{10}^{\otimes k})$.
We will examine it in detail in the next section.

\section{Centralizer algebra}
In the previous section,
we have seen that
the dimensions of ${\cal A}_k = \mbox{End}_{H_1}(V_{10}^{\otimes k})$
($k=1,2,\ldots$)
are 1, 2, 5, 15, 51, 187, 715,\ldots.
According to ``The On-Line Encyclopedia of Integer Sequences''\cite{Sl},
these terms coincide with the fist few terms of the expression
$(3\cdot 2^{k-2} + 2^{2k-3} + 1)/3$.
This is indeed the case.
In order to prove this,
we calculate the size of each simple component of ${\cal A}_k$.

Let $d^{(i)}_j$ be the multiplicity of $\rho_{i}$
in the tensor representation $\rho_{10}^j$,
which coincides with the size of the corresponding
simple component of ${\cal A}_j$.
By the Bratteli diagram of ${\cal A}_j$ given in the previous section,
we have the recursive formulae as follows.
First note that the irreducible representations
$\rho_8$, $\rho_{10}$, $\rho_{16}$ of $H_1$
again appear in the bottom of the diagram,
as well as the 5-th row of the diagram.
This implies that the diagram periodically grows up as $k$ increases.
The iteration is as follows:
\[
[\rho_8, \rho_{10}, \rho_{16}]\rightarrow
[\rho_4, \rho_3, \rho_{14}, \rho_{13}, \rho_6]\rightarrow
[\rho_9, \rho_7, \rho_{15}]\rightarrow
[\rho_1, \rho_2, \rho_{11}, \rho_{12}, \rho_5]\rightarrow\cdots.
\]
Hence based on the Bratteli diagram of the 9-th row from the 5-th row,
we can obtain the following recursive formulae:
\begin{align*}
&\left\{
	\begin{array}{rcl}
		d^{(4)}_{4\ell+2} &=& d^{(8)}_{4\ell+1},\\ 
		d^{(3)}_{4\ell+2} &=& d^{(10)}_{4\ell+1},\\
		d^{(14)}_{4\ell+2}
		&=& d^{(8)}_{4\ell+1}+d^{(16)}_{4\ell+1},\\
		d^{(13)}_{4\ell+2}
		&=& d^{(10)}_{4\ell+1}+d^{(16)}_{4\ell+1},\\
		d^{(6)}_{4\ell+2} &=& d^{(16)}_{4\ell+1},
	\end{array}
\right.\\
&\left\{
	\begin{array}{rcl}
		d^{(9)}_{4\ell+3} &=& d^{(4)}_{4\ell+2}+d^{(14)}_{4\ell+2},\\ 
		d^{(7)}_{4\ell+3} &=& d^{(3)}_{4\ell+2}+d^{(13)}_{4\ell+2},\\
		d^{(15)}_{4\ell+3}
		&=&
		d^{(14)}_{4\ell+2}+d^{(13)}_{4\ell+2}+d^{(6)}_{4\ell+2},\\
	\end{array}
\right.\\
&\left\{
	\begin{array}{rcl}
		d^{(1)}_{4(\ell+1)} &=& d^{(9)}_{4\ell+3},\\ 
		d^{(2)}_{4(\ell+1)} &=& d^{(7)}_{4\ell+3},\\
		d^{(11)}_{4(\ell+1)}
		&=& d^{(7)}_{4\ell+3}+d^{(15)}_{4\ell+3},\\
		d^{(12)}_{4(\ell+1)}
		&=& d^{(9)}_{4\ell+3}+d^{(15)}_{4\ell+3},\\
		d^{(5)}_{4(\ell+1)} &=& d^{(15)}_{4\ell+3},
	\end{array}
\right.\\
&\left\{
	\begin{array}{rcl}
		d^{(8)}_{4(\ell+1)+1}
		&=& d^{(2)}_{4(\ell+1)}+d^{(11)}_{4(\ell+1)},\\ 
		d^{(10)}_{4(\ell+1)+1}
		&=& d^{(1)}_{4(\ell+1)}+d^{(12)}_{4(\ell+1)},\\
		d^{(16)}_{4(\ell+1)+1}
		&=&
		d^{(11)}_{4(\ell+1)}+d^{(12)}_{4(\ell+1)}+d^{(5)}_{4(\ell+1)}.
	\end{array}
\right.
\end{align*}
Note also that if we allow the possibility $d_j^{(i)}=0$,
the recursions above are still valid even for the 1-st to 4-th row.
Hence we have the following:
\begin{align}\label{eq:rec8}\nonumber
d^{(8)}_{4\ell+1} &= d^{(2)}_{4\ell}+d^{(11)}_{4\ell}\\\nonumber
&= d^{(7)}_{4\ell-1}+(d^{(7)}_{4\ell-1}+d^{(15)}_{4\ell-1})\\\nonumber
&= 2d^{(7)}_{4\ell-1}+d^{(15)}_{4\ell-1}\\\nonumber
&= 2(d^{(3)}_{4\ell-2}+d^{(13)}_{4\ell-2})
+(d^{(14)}_{4\ell-2}+d^{(13)}_{4\ell-2}+d^{(6)}_{4\ell-2})\\\nonumber
&= 2d^{(3)}_{4\ell-2}+d^{(14)}_{4\ell-2}
+3d^{(13)}_{4\ell-2}+d^{(6)}_{4\ell-2}\\\nonumber
&= 2d^{(10)}_{4\ell-3}+(d^{(8)}_{4\ell-3}+d^{(16)}_{4\ell-3})
+3(d^{(10)}_{4\ell-3}+d^{(16)}_{4\ell-3})+d^{(16)}_{4\ell-3}\\
&= d^{(8)}_{4(\ell-1)+1}+5d^{(10)}_{4(\ell-1)+1}+5d^{(16)}_{4(\ell-1)+1}\ 
(\ell > 0).
\end{align}
Similarly we have
\begin{equation}\label{eq:rec10}
d^{(10)}_{4\ell+1}
= 5d^{(8)}_{4(\ell-1)+1}+d^{(10)}_{4(\ell-1)+1}+5d^{(16)}_{4(\ell-1)+1}
\end{equation}
and
\begin{equation}\label{eq:rec16}
d^{(16)}_{4\ell+1}
= 5d^{(8)}_{4(\ell-1)+1}+5d^{(10)}_{4(\ell-1)+1}+11d^{(16)}_{4(\ell-1)+1}
\end{equation}
From the recursion \eqref{eq:rec8}, \eqref{eq:rec10} and \eqref{eq:rec16},
and the initial condition $(d^{(8)}_1, d^{(10)}_1, d^{(16)}_1) = (0,1,0)$,
we obtain
\begin{align*}
d^{(8)}_{4\ell+1}
&=-\frac{(-4)^{\ell}}{2}+\frac{1}{3}+\frac{{16}^{\ell}}{6},\\
d^{(10)}_{4\ell+1}
&=\frac{(-4)^{\ell}}{2}+\frac{1}{3}+\frac{{16}^{\ell}}{6},\\
d^{(16)}_{4\ell+1}
&=-\frac{1}{3}+\frac{{16}^{\ell}}{3}.
\end{align*}
By the initial recursion formulae, we immediately obtain
\begin{align*}
d^{(4)}_{4\ell+2}
&=-\frac{(-4)^{\ell}}{2}+\frac{1}{3}+\frac{{16}^{\ell}}{6},\\
d^{(3)}_{4\ell+2}
&=\frac{(-4)^{\ell}}{2}+\frac{1}{3}+\frac{{16}^{\ell}}{6},\\
d^{(14)}_{4\ell+2}
&=-\frac{(-4)^{\ell}}{2}+\frac{{16}^{\ell}}{2},\\
d^{(13)}_{4\ell+2}
&=\frac{(-4)^{\ell}}{2}+\frac{{16}^{\ell}}{2},\\
d^{(6)}_{4\ell+2}
&=-\frac{1}{3}+\frac{{16}^{\ell}}{3},
\end{align*}
\begin{align*}
d^{(9)}_{4\ell+3}
&=-(-4)^{\ell}+\frac{1}{3}+\frac{2\cdot{16}^{\ell}}{3},\\
d^{(7)}_{4\ell+3}
&=(-4)^{\ell}+\frac{1}{3}+\frac{2\cdot{16}^{\ell}}{3},\\
d^{(15)}_{4\ell+3}
&=-\frac{1}{3}+\frac{4\cdot{16}^{\ell}}{3}
\end{align*}
and
\begin{align*}
d^{(1)}_{4(\ell+1)}
&=-(-4)^{\ell}+\frac{1}{3}+\frac{2\cdot{16}^{\ell}}{3},\\
d^{(2)}_{4(\ell+1)}
&=(-4)^{\ell}+\frac{1}{3}+\frac{2\cdot{16}^{\ell}}{3},\\
d^{(11)}_{4(\ell+1)}
&=(-4)^{\ell}+2\cdot{16}^{\ell},\\
d^{(12)}_{4(\ell+1)}
&=-(-4)^{\ell}+2\cdot{16}^{\ell},\\
d^{(5)}_{4(\ell+1)}
&=-\frac{1}{3}+\frac{4\cdot{16}^{\ell}}{3}.
\end{align*}
Thus we have obtained the size of each simple component of ${\cal A}_k$.
If we apply simple considerations to the order of simple components,
the sizes are uniformly described as follows.
\begin{thm}
Let ${\cal A}_k = \End_{H_1}(V_{10}^{\otimes k})$
be a centralizer algebra of $H_1$ in $V_{10}^{\otimes k}$,
where $H_1$ acts on $V_{10}$ diagonally.
Then ${\cal A}_k$ has the following multi-matrix structure.
\[
	{\cal A}_{k}
	\cong
	\left\{
	\begin{array}{l}
		M_{d_+(k)}(\Comp)
		\oplus M_{d_-(k)}(\Comp)
		\oplus M_{d_0(k)}(\Comp)
		\\
		\quad\mbox{if $k = 2m-1$},\\
		M_{d_+(k)}(\Comp)
		\oplus M_{d_-(k)}(\Comp)
		\oplus M_{d_0(k)}(\Comp)
		\oplus M_{e_+(k)}(\Comp)
		\oplus M_{e_-(k)}(\Comp)
		\\
		\quad\mbox{if $k = 2m$},
	\end{array}
	\right.
\]
where
\begin{align*}
	d_{\pm}(k) &= \pm 2^{m-2}+\frac{1}{3}+\frac{2\cdot 4^{m-2}}{3},\\
	d_0(k)& = -\frac{1}{3} + \frac{4^{m-1}}{3}
\end{align*}
and
\[
	e_{\pm} =\pm 2^{m-2}+2\cdot 4^{m-2}.
\]
\end{thm}
Calculating the square sum of the dimensions of
the simple components of ${\cal A}_k$ in cases $k=2m-1$ and $k=2m$,
we finally obtain the following corollary as we expected.
\begin{cor}
\[
\dim {\cal A}_k = 2^{k-2} + \frac{2^{2k-3}}{3} + \frac{1}{3}.
\]
\end{cor}

Again by the web cite~\cite{Sl},
this result suggests that
the basis of the centralizer algebras of $H_1$
could be described in terms of
the symmetric polynomials in 4 noncommuting variables~\cite{Be}
and/or the universal embedding
of the symplectic dual polar space $DS_p(2k,2)$~\cite{Bl}.
It would be interesting that these points become clear.

\bigskip

{\bf Acknowledgment}. 
This work was supported by JSPS KAKENHI Grant Number 25400014.

{\sc Department of Mathematical Sciences, University of the Ryukyus,
Okinawa,
903-0213, JAPAN}

{\it E-mail address}: kosuda@math.u-ryukyu.ac.jp

\bigskip

{\sc Graduate School of Natural Science and Technology,
Kanazawa University,
Ishikawa, 920-1192
Japan}

{\it E-mail address}: oura@se.kanazawa-u.ac.jp

\end{document}